\documentclass[]{interact}

\usepackage{epstopdf}
\usepackage[caption=false]{subfig}
\usepackage{dsfont}
\usepackage[numbers,sort&compress]{natbib}
\bibpunct[, ]{[}{]}{,}{n}{,}{,}

\usepackage{natbib}

\DeclareMathOperator{\Var}{Var}

\renewcommand\bibfont{\fontsize{10}{12}\selectfont}

\theoremstyle{plain}
\newtheorem{theorem}{Theorem}[section]
\newtheorem{lemma}[theorem]{Lemma}

\theoremstyle{definition}

\theoremstyle{remark}
\newtheorem{remark}{Remark}

\begin{document}


\title{The Bahadur representation of sample quantiles for associated sequences}

\author{
	\name{Lahcen DOUGE \thanks{lahcen.douge@uca.ac.ma} }
	\affil{FSTG, University Cadi Ayyad, B.P. 549 Marrakech, Morocco }
}

\maketitle

\begin{abstract}
	In this paper, the Bahadur representation of sample quantiles based on associated sequences is established under  polynomially decaying of covariances. The rate of approximation depends on the covariances decay degree and  becomes close to the optimal rate obtained  under independence when the covariances decrease fastly to $0$.
\end{abstract}


\section{Introduction}

Let $(X_n)_{n\geq1}$ be a sequence of strictly stationary associated random variables. Assume that $X_1$ has continuous  distribution function $F$ and let $f$ and $Q$ denote the associated density and quantile function, respectively. For $0<p<1$, denote by $ \xi_p=Q(p)$ the $p$th quantile of $F$. Given a sample $X_1,\ldots, X_n$,  define the  empirical distribution function 
\begin{eqnarray*}
F_n(x)= \frac{1}{n}\sum_{i=1}^{n} \mathds{1}_{\{X_{i}\leq x\}},\quad x\in\mathbb{R},
\end{eqnarray*}
where $\mathds{1}_{A}$ denotes the indicator function of a set $A$ and let $\xi_{n,p}=\inf\{x : F_n(x)\geq p\}$ the $p$th sample quantile. With $U_i=F(X_i)$, $i\geq 1$, define, for each $0\leq t \leq 1$, the uniform empirical distribution
\begin{eqnarray*}
E_n(t)= \frac{1}{n}\sum_{i=1}^{n} \mathds{1}_{\{U_{i}\leq t\}}=F_n\left(Q(t)\right).
\end{eqnarray*}


Bahadur representation is  useful to establish consistency and asymptotic normality results for  sample quantiles.
\citet{bahadur} has established the asymptotic representation for sample quantile via the  empirical distribution function based on independent and identically distributed (i.i.d.) random variables. \citet{kiefer} provided exact rates in the Bahadur representation for i.i.d. sequences. The extensions of Bahadur type representations  by relaxing the assumption
of independence have been studied by a number of
authors. Especially, \citet{sen} obtained  similar results to Bahadur's
one for stationary $\phi$-mixing processes.  \citet{yoshihara} provided some generalizations of Sen's results 
for $\phi$-mixing and $\alpha$-mixing sequences. 
 Recent works on Bahadur representation for sample quantiles, among many others,  include, for example, \citet{zhang2014}, \citet{wang2016}, \citet{xing2019} and \citet{wu2021} for mixing sequences, and \citet{xu2013} for negatively associated sequences. Other extensions of Bahadur' representations  under weak dependence were investigated in \citet{wu}, \citet{sun2006}, \citet{yang2014}, \citet{yang2019}, \citet{kong2017} and \citet{wu2021nonparametric}.

The object of the present paper is to show the Bahadur representation of sample quantiles for associated sequences under the polynomial decay of the covariances. The rate of approximation obtained is close to the optimal one when the covariance decay degree is large enough.

Companies in the finance and insurance sector manage significant risks of various kinds. One of the tools involved in this management is the measurement of risk 
by using the quantile function or what is known as  Value-at-risk (VaR). 
Usually, the risks of the insurance and finance portfolio are assumed independent. However, there are practical situations for which this assumption is not appropriate. Some examples of associated risks in insurance  were given in \citet{cossette} and \citet{denuit2006}.

Before stating our results, we recall the notion of  association, which was introduced by \citet{esary1967}. A sequence $(X_n)_{n\geq1}$ of real-valued random variables is said to be associated if, for every finite subcollection $X_{i_1},\ldots,X_{i_n}$ and every pair of coordinatewise nondecreasing functions $f_1$, $f_2$ : $\mathbb{R}^{n}\rightarrow   \mathbb{R}$,
\begin{eqnarray*}
	\mathrm{Cov}\big( f_1(X_{i_1},\ldots,X_{i_n}), f_2(X_{i_1},\ldots,X_{i_n})\big) \geq 0,
\end{eqnarray*}
whenever the covariance is defined. The main advantage of dealing with associated random variables is that the conditions of limit theorems are based on covariance structure which is easier to evaluate than mixing coefficients. 
 The fundamental results in this research domain with extensive bibliographical references can be found
in the monograph of \citet{bulinski2007}. 
%
%
%

\section{Main results}
To prove the  main results, we need the following assumptions.
\begin{enumerate}
	\item[(A.1)] 
	 $f$ possesses a bounded derivative $f^{'}$ in a neighborhood of $\xi_p$.
	\item[(A.2)] There exist constants $b_0\geq 0$ and $ b >0$ such that for all $k\geq 1$,
	$$\mathrm{Cov}\left( X_1, X_{k+1}\right) \leq b_0 k^{-b}.$$
\end{enumerate}
Throughout this paper, denote by $c$, $c_1$ and $c_2$  strictly positive constants  whose values are allowed to change in each appearance and consider that a sequence of random variables $(Z_n)_{n\geq 1}$ is said to be $O_{a.s.}(r_n)$ if $\frac{Z_n}{r_n}$ is almost surely bounded.
\begin{theorem}\label{thm1} Suppose that  (A.2) holds for  $b>3$. If  $f$ is bounded and $f(\xi_{p})>0$, then 
	\begin{eqnarray}\label{thm1eq1}
	\left| \xi_{n,p}-\xi_p\right| =O_{a.s.}\left( n^{-\frac{1}{2}}\log^{\delta} n\right),
	\end{eqnarray}
	where $\delta >\frac{3}{2b}$.
\end{theorem}

\begin{theorem}\label{thm2}  Suppose that  (A.2) holds for $b>\frac{5+\sqrt{17}}{2}$. If  $f$ is bounded, then
	\begin{eqnarray}\label{th2eq1}
	\sup_{t\in J_n} \Big|\big (E_n(t)- t\big)- \big(E_n(p)-p)\big)\Big| = O_{a.s.}\left(n^{-\frac{1}{2}-\frac{\beta_b }{4}} \log^{\gamma}n\right),
	\end{eqnarray}
	where $J_n=\left\{t : |t - p|\leq c\, n^{-\frac{1}{2}}\log^{\delta} n\right\}$, $\beta_b=\frac{b-3}{b-1}$ and  $\gamma> \frac{3 \beta_b }{4 b}+\frac{1}{b}$.
\end{theorem}

\begin{theorem}\label{thm3} Suppose that  (A.1) and (A.2) hold for some $b>\frac{5+\sqrt{17}}{2}$.  If  $f$ is bounded and $f(\xi_p)>0$, then
	\begin{eqnarray}\label{thm2eq1}
	\xi_{n,p}-\xi_p=\frac{p-F_n(\xi_p)}{f(\xi_p)}	 + O_{a.s.}\left(n^{-\frac{1}{2}-\frac{\beta_b }{4}} \log^{\gamma}n\right).
	\end{eqnarray}
\end{theorem}

\begin{remark}
	We remark that the rate  in (\ref{thm2eq1}) becomes  close to the  optimal bound $n^{-\frac{3}{4}}$ when $b$ is large enough. 
	
\end{remark}

%
\begin{remark}
	 Following the steps of the proof of Theorem 2.4 in \citet{xing2019}, we can  deduce from Theorem \ref{thm3} the uniformly asymptotic normality of sample quantiles for associated random variables. However the best rate derived will be of the order of $n^{-\frac{\beta_b }{4}} \log^{\gamma}n$ which is slower  than the rate  obtained in \citet{douge2020}.
\end{remark}

\section{Proofs of  main results}



\textbf{Proof of Theorem \ref{thm1}.} Let $a_n=c \,  n^{-\frac{1}{2}}\log^{\delta} n$ and consider, for any $k\geq 0$,  the event
\begin{eqnarray*}
	A_k=\bigcup_{2^k\leq n <2^{k+1}}\Big(| \xi_{n,p}-\xi_{p}| > a_n\Big).
\end{eqnarray*}
Clearly
\begin{eqnarray*}
	\big(| \xi_{n,p}-\xi_{p}| > a_n\big)&=& \big(\xi_{n,p}>\xi_{p}+a_n\big)\cup \big(\xi_{n,p}<\xi_{p}-a_n\big).
\end{eqnarray*}
Now
\begin{eqnarray*}\label{eq4}
	\big(\xi_{n,p}>\xi_{p}+a_n\big) &=&\left(\sum_{i=1}^{n} \mathds{1}_{\{X_i > \xi_p+a_n\}}> n(1-p) \right)\\
	&=&\left(\sum_{i=1}^{n} \big[ \mathds{1}_{\{X_i > \xi_p+a_n\}}- E\mathds{1}_{\{X_i > \xi_p+a_n\}}\big]> n \delta_1 \right)
\end{eqnarray*}
and
\begin{eqnarray*}\label{eq4}
	\big(\xi_{n,p}\leq \xi_{p}-a_n\big) =\left(\sum_{i=1}^{n}  \big[\mathds{1}_{\{X_i \leq \xi_p-a_n\}}- E\mathds{1}_{\{X_i \leq \xi_p-a_n\}}\big]\geq n \delta_2\right),
\end{eqnarray*}
where $\delta_1=F(\xi_{p}+a_n)-p$ and $\delta_2=p-F(\xi_{p}-a_n)$. Since $F$ is continuous at $\xi_p$, $F(\xi_p)=p$. For $n$ large enough, by  Taylor's expansion formula there exists $\eta_1\in \left]0,  a_n\right[$ such that
\begin{eqnarray*}
	F(\xi_p+a_n)-p&=&f(\xi_p)a_n+\frac{1}{2}\, f^{'}(\xi_p + \eta_1)a_n^2.
\end{eqnarray*}
For $n$ large enough, by (A.1) we have
\begin{eqnarray*}
	\sup_{|x| \leq a_n}\big|f^{'}(\xi_p+x)\big| \leq c_1<\infty
\end{eqnarray*}
and from this  we deduce that
\begin{eqnarray*}\label{thm1eq22b}
c_2 \, a_n	\leq a_n  \big[ f(\xi_p)- \frac{1}{2} c_1 \, a_n\big]\leq \delta_1.
\end{eqnarray*}
It follows that
\begin{eqnarray*}\label{thm1eq22b}
	\bigcup_{2^k\leq n <2^{k+1}} \big(\xi_{n,p}>\xi_{p}+a_n\big)&\subset&  \bigcup_{2^k\leq n <2^{k+1}} \left\{\sum_{i=1}^{n} V_i > c \, n^{\frac{1}{2}} \log^{\delta} n\right\}\\
	&\subset&  \left\{\max_{2^k\leq n \leq 2^{k+1}}\left|\sum_{i=1}^{n} V_i\right|> c \,  2^{\frac{k}{2}} k^{\delta} \right\},
\end{eqnarray*}
where $V_i=\mathds{1}_{\{X_i > \xi_p+a_n\}}- E\mathds{1}_{\{X_i > \xi_p+a_n\}}$, $i\geq 1$. We will use now the maximal inequality for associated random variables in Lemma \ref{lemA3}. First, observe that by the stationarity of the associated sequence $(V_i)_{i\geq 1}$, for any $m\geq 2$

\begin{eqnarray}\label{Pthm1}
s_m^2 := \Var\left(\sum_{i=1}^{m} V_i \right)&=& m \Var (V_1)+2\sum_{l=2}^{m} (m-l+1) \mathrm{Cov}(V_1,V_l)\nonumber\\
&\leq & m+2 m \sum_{l=2}^{\infty} \mathrm{Cov}(V_1,V_l).
\end{eqnarray}
and
\begin{eqnarray}\label{Pthm3}
m \Var (V_1)\leq s_m^2 .
\end{eqnarray}
We apply now Lemma \ref{lemA4}  and Lemma \ref{lemA6} for the associated sequence $(U_i)_{i\geq 1}$ to   get 
\begin{eqnarray}\label{Pthm2}
\mathrm{Cov}(V_1,V_l)&=& \mathrm{Cov}\left(\mathds{1}_{\{U_1> F(\xi_p+a_n)\}},\mathds{1}_{\{U_l> F(\xi_p+a_n)\}}\right)\nonumber\\
&\leq & 4\,  \mathrm{Cov}(U_1,U_l)^{\frac{1}{3}}\nonumber\\
&\leq& 4 \|f\|_{\infty}^{\frac{2}{3}} \mathrm{Cov}(X_1,X_l)^{\frac{1}{3}}.
\end{eqnarray}
By condition (A.2), with $b>3$,  we conclude from(\ref{Pthm1})-(\ref{Pthm2}) that,  for $n$ large enough,
\begin{eqnarray}\label{Pthm4}
c_1m\leq s_m^2\leq c_2 m.
\end{eqnarray}
 Since the associated sequence $(V_i)_{i\geq 1}$ is  centered, by Lemma \ref{lemA3}  and (\ref{Pthm4}), it follows that, for $n$ large enough, 
\begin{eqnarray*}\label{thm1eq22b}
	P \left(\max_{2^k\leq n \leq 2^{k+1}}\left|\sum_{i=1}^{n} V_i \right|> c \, 2^{\frac{k}{2}} k^{\delta} \right) 
		&\leq & P \left(\max_{1\leq j \leq 2^{k+1}}\left|\sum_{i=1}^{j} V_i \right|> c \, k^{\delta}  s_{2^{k+1}} \right)\\
	&\leq & P \left(\left|\sum_{i=1}^{2^{k+1}} V_i\right|> c \, k^{\delta}  s_{2^{k+1}} \right).
\end{eqnarray*}
By Markov inequality, (\ref{Pthm2}), (\ref{Pthm4}) and Lemma \ref{lemA2}, we get
\begin{eqnarray*}\label{thm1eq22b}
	P \left(\left|\sum_{i=1}^{2^{k+1}} V_i\right|> c \, k^{\delta}  s_{2^{k+1}} \right)&\leq & P \left(\left|\sum_{i=1}^{2^{k+1}} V_i\right|> c \, k^{\delta}  2^{\frac{k}{2}}\right)\\
	&\leq & c \, 2^{-\frac{bk}{3}} k^{-\frac{2b\delta}{3}} E\Bigg|\sum_{i=1}^{2^{k+1}} V_i\Bigg|^{\frac{2b}{3}}\\
	&\leq & c \,   k^{-\frac{2b\delta}{3}}.
\end{eqnarray*}
We thus obtain
\begin{eqnarray*}\label{thm1eq22b}
	P \Big(	\bigcup_{2^k\leq n <2^{k+1}} \big(\xi_{n,p}>\xi_{p}+ a_n\big)\Big) \leq  c \,k^{-\frac{2b\delta}{3}}.
\end{eqnarray*}
Similarly,
\begin{eqnarray*}\label{thm1eq22b}
	P \Big(	\bigcup_{2^k\leq n <2^{k+1}} \big(\xi_{n,p}<\xi_{p}-a_n\big)\Big) \leq  c \, k^{-\frac{2b\delta}{3}}.
\end{eqnarray*}
Therefore, as $\delta>\frac{3}{2b}$,  we have $\sum_{k=1}^{\infty} P(A_k)<\infty$, which completes the proof of the theorem by applying Borel Cantelli's lemma.\\

\noindent To prove Theorem \ref{thm2}, we need the following lemma.
\begin{lemma}  \label{lem1}
	Let $q>4$. Suppose that $f$ is bounded and (A2) holds for some $b>q-1$. If  $ t-s > 2 n^{\alpha_q}$, where $\alpha_q= \frac{(-q+4+2\eta)(q-1)}{(q+2)(q-3)}<0$ for some $\eta>0$, then  		
	\begin{eqnarray*}
		E\Big|E_n(t)-E_n(s)-(t-s)\Big|^{q}\leq c\, n^{-\frac{q}{2}} (t-s)^{\frac{q(q-3)}{2(q-1)}}.
	\end{eqnarray*}
\end{lemma}
\textbf{Proof.} By Lemma \ref{lemA6}, we have
\begin{eqnarray}\label{lem1eq1}
\mathrm{Cov}(U_1, U_n)\leq \|f\|_{\infty}^{2}\mathrm{Cov}(X_1, X_n)=O(n^{-b}).
\end{eqnarray}
Set $Z_i=\mathds{1}_{\{s<U_i\leq t\}}-(t-s)$, $i\geq 1$.  $(U_i)_{i\geq 1}$ is a stationary sequence of associated  uniform $[0,1]$ random variables. Observe that, for any $i\geq 1$, $|E Z_1 Z_i| \leq t-s$. By Lemma \ref{lemA4}, (A2) and (\ref{lem1eq1}), we get
\begin{eqnarray*}
	\sum_{i=1}^{n}|E Z_1 Z_i| &\leq & (t-s)^{1-\frac{3}{q-1}} \sum_{i=1}^{\infty}|E Z_1 Z_i|^{\frac{3}{q-1}}\\
	&\leq & 4 \,(t-s)^{1-\frac{3}{q-1}} \sum_{i=1}^{\infty} (t-s)^{\frac{1}{q-1}} \mathrm{Cov}(U_1, U_i)^{\frac{1}{q-1}}\\
	&\leq & c\, (t-s)^{1-\frac{2}{q-1}}.
\end{eqnarray*}
Now, by Lemma \ref{lemA5}, we obtain, for $t-s>2 n^{\frac{-q+1+\eta}{q+2}}$,
\begin{eqnarray*}
	E\Big|E_n(t)-E_n(s)-(t-s)\Big|^{q}\leq \frac{c}{n^{q}} \Big(n^{\frac{q(3+\eta)}{q+2}}+\big(n (t-s)^{\frac{q-3}{q-1}}\big)^{\frac{q}{2}}\Big).
\end{eqnarray*}
To complete the proof of the lemma, it suffices to check first that  if  $t-s  > 2 n^{\alpha_q}$  we have
$n^{\frac{q(3+\eta)}{q+2}} < c \big(n (t-s)^{\frac{q-3}{q-1}}\big)^{\frac{q}{2}}$, and second that $\alpha_q>\frac{-q+1+\eta}{q+2}$ for a suitable choice of $\eta$.\\\\

\noindent\textbf{Proof of Theorem \ref{thm2}.} For any $n\geq 1$, set $Y_n(t)=E_n(t)- t$, $0\leq t \leq 1$. For any $k\geq 1$, define the event
\begin{eqnarray*}\label{thm2eq0}
	B_k =\bigcup_{2^k\leq n < 2^{k+1}}\Big\{	\sup_{t\in J_n} \big|Y_n(t)- Y_n(p)\big| \geq c\, n^{-\frac{1}{2}- \frac{\beta_b}{4}}\log^{\gamma} n \Big\}.
\end{eqnarray*}
 Let $\theta_n=c_1\, n^{-\frac{1}{2}- \frac{\beta_b}{4}}$
and $\gamma_n =\left[c_2\frac{n^{-\frac{1}{2}}\log^{\delta} n}{\theta_n}\right]+1$, where $[x]$ is the integer part of $x$. Now, for each $2^k\leq n \leq 2^{k+1}$, we have
\begin{eqnarray*}\label{thm2eq0}
	\sup_{t\in J_n} \big|Y_n(t)- Y_n(p)\big| &\leq&  \sup_{p -\theta_n \gamma_n\leq t \leq p+ \theta_n \gamma_n} \big|Y_n(t)- Y_n(p)\big| \\
	&=&\max_{1\leq j \leq \gamma_n} \sup_{p +(j-1)\theta_n\leq t \leq p +j \theta_n } \big|Y_n(t)- Y_n(p)\big|\\
	&&\,\, \vee 
	\max_{1\leq j \leq \gamma_n} \sup_{p -j \theta_n \leq t\leq p -(j-1)\theta_n } \big|Y_n(t)- Y_n(p)\big| \\
	&=:& \Gamma_1 \vee \Gamma_2,
\end{eqnarray*}
where $x\vee y=max\{x,y\}$. Since $E_n$ is increasing, we get 
\begin{eqnarray*}\label{thm2eq0}
	\Gamma_1 &\leq &\max_{1\leq j \leq \gamma_n}  \Big|\big(E_n(p+j \theta_n)-(p+(j-1) \theta_n)\big)- Y_n(p)\Big| \\
	&&\,\, \vee  \max_{1\leq j \leq \gamma_n}  \Big|\big(E_n(p+(j-1)\theta_n)- (p+j \theta_n)\big)- Y_n(p)\Big| \\
	&\leq& \Big(\max_{1\leq j \leq \gamma_n}  \big|Y_n(p+j \theta_n)- Y_n(p)\big|+\theta_n\Big) \\
	&&\,\, \vee  \Big(\max_{1\leq j \leq \gamma_n}  \big|Y_n(p+(j-1)\theta_n)- Y_n(p)\big|+\theta_n\Big) \\
	&=& \max_{1\leq j \leq \gamma_n}  \big|Y_n(p+j \theta_n)- Y_n(p)\big|+\theta_n.
\end{eqnarray*}
Likewise,
\begin{eqnarray*}\label{thm2eq0}
	\Gamma_2 \leq  \max_{1\leq j \leq \gamma_n}  \big|Y_n(p-j \theta_n)- Y_n(p)\big|+\theta_n.
\end{eqnarray*}
Then we have
\begin{eqnarray*}\label{thm2eq0}
	B_k \subset B_{k1} \cup B_{k2},
\end{eqnarray*}
where
\begin{eqnarray*}\label{thm2eq0}
	B_{k1} =\bigcup_{2^k\leq n < 2^{k+1}}\Big\{	\max_{1\leq j \leq \gamma_n}  \big|Y_n(p+j \theta_n)- Y_n(p)\big|\geq c\, n^{-\frac{1}{2}- \frac{\beta_b}{4}}\log^{\gamma} n \Big\}
\end{eqnarray*}
and 
\begin{eqnarray*}\label{thm2eq0}
	B_{k2} =\bigcup_{2^k\leq n < 2^{k+1}}\Big\{	\max_{1\leq j \leq \gamma_n}  \big|Y_n(p - j \theta_n)- Y_n(p)\big|\geq c\, n^{-\frac{1}{2}- \frac{\beta_b}{4}}\log^{\gamma} n \Big\}.
\end{eqnarray*}
Now, for $0\leq s<t\leq 1$, denote
\begin{eqnarray*}\label{thm2eq0}
	\eta_i(s,t)= \big(\mathds{1}_{\{U_i\leq t\}}- t\big) - \big(\mathds{1}_{\{U_i\leq s\}}- s\big).
\end{eqnarray*}
Clearly
\begin{eqnarray*}\label{thm2eq0}
	\big|Y_n(p + j \theta_n)- Y_n(p)\big|= \frac{1}{n}\left|\sum_{i=1}^{n} \eta_i(p, p + j \theta_n)\right|= \frac{1}{n}\left|\sum_{i=1}^{n} \sum_{l=1}^{j}\eta_{i,l}\right|
\end{eqnarray*}
and 
\begin{eqnarray*}\label{thm2eq0}
	\big|Y_n(p - j \theta_n)- Y_n(p)\big|= \frac{1}{n}\left|\sum_{i=1}^{n} \eta_i(p-j \theta_n, p )\right|= \frac{1}{n}\left|\sum_{i=1}^{n} \sum_{l=1}^{j}\delta_{i,l}\right|,
\end{eqnarray*}
where
\begin{eqnarray*}\label{thm2eq0}
	\eta_{i,l}= \eta_i\big(p + (l-1)\theta_n, p + l \theta_n\big)\quad \text{and}\quad \delta_{i,l}= \eta_i\big(p - l \theta_n, p -  (l-1) \theta_n\big).
\end{eqnarray*}
 For an appropriate choice of $c_1$, we have, for every $n\geq 1$, $\theta_n > 2 n^{ \alpha_b }$, $\alpha_b= \frac{(-b+4+2\eta)}{(b+2)\beta_b}$. Since the sequence $(U_i)_{i\geq 1}$ is  strictly stationary. For every integers $n\geq 1$, $j\geq 1$ and $b_1, b_2 \in \mathbb{N}$, by Lemma   \ref{lem1}, for some $\tau>0$, we get
\begin{eqnarray*}\label{thm2eq0}
	E\left|\sum_{i=b_1+1}^{b_1+n} \sum_{l=b_2+1}^{j+b_2}\eta_{i,l}\right|^b &=& E\left|\sum_{i= 1}^{n} \sum_{l= b_2 +1}^{j + b_2}\eta_{i,l}\right|^b\\
	&=&   n^b E\big|E_n(p+ (j+b_2)\theta_n)-E_n(p + b_2\theta_n )- j\theta_n\big|^b\\
	&\leq& c\,  n^{\frac{b}{2}} (j \theta_n)^{\frac{b\beta_b }{2}}\\
	&\leq & \big(c \, n^{\frac{b}{2+\tau}  - \frac{b \beta_b}{2+\tau}(\frac{1}{2}+ \frac{\beta_b}{4})}j^{\frac{b\beta_b}{2+\tau}}\big)^{\frac{2+\tau}{2}}:=\big(g\big(R_{(b_1,b_2), (n,j)}\big)\big)^{\frac{2+\tau}{2}},
	\end{eqnarray*}
where	$R_{(b_1,b_2), (n,j)}=\big\{(i_1, i_2)\in \mathbb{N}^{2}: b_1+1\leq i_1 \leq b_1+n,\, b_2+1\leq i_2 \leq b_2+j\big\}$. According to $\frac{b}{2+\tau}  - \frac{b \beta_b}{2+\tau}(\frac{1}{2}+ \frac{\beta_b}{4})\geq 1$ and $\frac{b\beta_b}{2+\tau}\geq 1$, for $b>\frac{5+\sqrt{17}}{2}$ and a suitable choice of $\tau$, we conclude that $g$ is superadditive and hence by Lemma \ref{lemA1} it follows that
\begin{eqnarray*}\label{thm2eq0}
	E\left(\max_{1\leq n \leq 2^{k+1}}\max_{1\leq j \leq \gamma_{2^{k+1}}}\left|\sum_{i=1}^{n} \sum_{l=1}^{j}\eta_{i,l}\right|\right)^b &\leq&  \big(g\big(R_{(0,0), (2^{k+1}, \gamma_{2^{k+1}})}\big)\big)^{\frac{2+\tau}{2}}\\
	& = & c\Big(2^{k(\frac{1}{2}-\frac{\beta_b }{4})}k^{\frac{\delta \beta_b }{2}}\Big)^b.
\end{eqnarray*}
From this, for any given $\varepsilon>0$, we obtain
\begin{eqnarray*}\label{thm2eq0}
P(B_{k1}) &\leq &	P\left(\max_{1\leq n \leq 2^{k+1}}\max_{1\leq j \leq \gamma_{2^{k+1}}}\left|\sum_{i=1}^{n} \sum_{l=1}^{\gamma_{2^{k+1}}}\eta_{i,l}\right|\geq  \varepsilon 2^{k(\frac{1}{2}-\frac{\beta_b }{4})} k^{\gamma}\right)\\
	&\leq & c\,  2^{-kb(\frac{1}{2}-\frac{\beta_p }{4})} k^{-\gamma b}E\left(\max_{1\leq n \leq 2^{k+1}}\max_{1\leq j \leq \gamma_n}\left|\sum_{i=1}^{n} \sum_{l=1}^{\gamma_{2^{k+1}}}\eta_{i,l}\right|\right)^b\\
	&\leq & c\, k^{-\gamma b + \frac{\delta b\beta_b }{2}}.
\end{eqnarray*}
Similar arguments show  that
\begin{eqnarray*}\label{thm2eq0}
	P(B_{k2}) 	\leq c\, k^{-\gamma b + \frac{\delta b\beta_b }{2}}.
\end{eqnarray*}
Finally, since $\gamma> \frac{\delta \beta_b }{2}+\frac{1}{b}$, with $\delta$ is taken to be close enough  to
$\frac{3}{2 b}$, we have $\sum_{k=1}^{\infty} P(B_k) <\infty$ and (\ref{th2eq1}) follows by using the Borel-Cantelli lemma.\\

\noindent\textbf{Proof of Theorem \ref{thm3}.} Let $\tilde{\xi}_{n,p}$ denote the pth quantile of  $E_n$. We apply now Theorem \ref{thm1} to the associated sequence $(U_i)_{i\geq 1}$ and we get
\begin{eqnarray*}\label{thm1eq1}
	\left| \tilde{\xi}_{n,p}- p\right| = O_{a.s.}\left(n^{-1/2}\log^{\delta} n\right).
\end{eqnarray*}
From (\ref{th2eq1}) we see that
\begin{eqnarray*}
	\left|\big (E_n(\tilde{\xi}_{n,p})- \tilde{\xi}_{n,p}\big)- \big(E_n(p)-p)\big)\right| = O_{a.s.}\left(n^{-\frac{1}{2}-\frac{\beta_b }{4}} \log^{\gamma}n\right).
\end{eqnarray*}
On the other hand, on noting  that (see \citet*{sen})
\begin{eqnarray*}
	E_n (\tilde{\xi}_{n,p})=\frac{r}{n}= p+O\left(\frac{1}{n}\right),\quad r=[np]+1,
\end{eqnarray*}
we obtain
\begin{eqnarray*}
	\Big|\big (p- \tilde{\xi}_{n,p}\big)- \big(E_n(p)-p)\big)\Big| = O_{a.s.}\left(n^{-\frac{1}{2}-\frac{\beta_b }{4}} \log^{\gamma}n\right).
\end{eqnarray*}
Since $F$ is continuous and increasing, we can check easily that $E_n(p)=F_n(\xi_{p})$ and  $\tilde{\xi}_{n,p}=F(\xi_{n,p})$. Consequently
\begin{eqnarray}\label{pthm3eq1}
\Big|\big (p- F(\xi_{n,p})\big)- \big(F_n(\xi_{p})-p)\big)\Big| = O_{a.s.}\left(n^{-\frac{1}{2}-\frac{\beta_b }{4}} \log^{\gamma}n\right).
\end{eqnarray}
Now by Taylor's expansion, we get 
\begin{eqnarray*}
	F(\xi_{n,p})=p+f(\xi_p)(\xi_{n,p}-\xi_p) +\frac{1}{2} f^{'}\big(\xi_p+\theta (\xi_{n,p}-\xi_p)\big)(\xi_{n,p}-\xi_p)^2,
\end{eqnarray*}
where $|\theta|<1$. By Theorem \ref{thm1} and (A1),
\begin{eqnarray*}
	F(\xi_{n,p})- p - f(\xi_p)(\xi_{n,p}-\xi_p) = O_{a.s.}\left(n^{-1} \log^{2\delta}n\right),
\end{eqnarray*}
which, together with (\ref{pthm3eq1}), leads to
\begin{eqnarray*}
	\Big|\big (f(\xi_p)(\xi_{n,p}-\xi_p)+ \big(F_n(\xi_{p})-p)\big)\Big| = O_{a.s.}\left(n^{-\frac{1}{2}-\frac{\beta_b }{4}} \log^{\gamma}n\right).
\end{eqnarray*}
The proof is completed.


\appendix
	\section{}
\renewcommand{\thelemma}{\Alph{section}\arabic{lemma}}

\begin{lemma}[ \citet{birkel1988b}, Lemma 3.1]\label{lemA6}	
	Let $A$ and $B$ be finite sets and let $(X_j)_{j\in A\cup B}$ be associated random variables. Then for all  real-valued partially differentiable functions $h_1$, $h_2$ with bounded partial derivatives, there holds
	\begin{eqnarray*}
		\Big| \mathrm{Cov}\big( h_1\big( (X_i)_{i\in A}\big) , h_2\big( (X_j)_{j\in B}\big) \big)  \Big|  \leq \sum_{i\in A}\sum_{j\in B} \left\|\frac{ \partial h_1}{\partial x_i} \right\| _{\infty}\left\| \frac{ \partial h_2}{\partial x_j}\right\| _{\infty}  \mathrm{Cov}(X_i,Y_j).
	\end{eqnarray*}
\end{lemma}

\begin{lemma}[ \citet{shao}, Lemma 5.1] \label{lemA4} 	Let $X$ and $Y$ be associated random variables with a common uniform distribution over $[0,1]$. Then for any $0\leq s < t \leq 1$,
	\begin{eqnarray*}
		\Big|\mathrm{Cov}\big(\mathds{1}_{\{s<X\leq t\}}, \mathds{1}_{\{s<Y\leq t\}}\big)\Big| \leq 4 (t-s)^{\frac{1}{3}} \big(\mathrm{Cov}(X,Y)\big)^{\frac{1}{3}}.
	\end{eqnarray*}
\end{lemma}

\begin{lemma}[ \citet{shao},  (5.27)] \label{lemA5} 	Let $q>2$. Let $(U_i)_{i\geq 1}$ be a stationary associated sequence of uniform $[0,1]$ random variables and let  $Z_i=\mathds{1}_{\{s<U_i\leq t\}}-(t-s)$, $0\leq s<t\leq 1$. If $t-s > 2 n^{\frac{-q+1+\eta}{q+2}}$ and
	$\mathrm{Cov}(U_1,U_n)=O(n^{-b})$ for some $ b>q-1$, then, for any $\eta>0$,  there exists some positive constant $K_{\eta}$ for which
	\begin{eqnarray*}
		E\Big|\sum_{i=1}^{n} Z_i\Big|^{q}\leq K_{\eta} \Bigg\{n^{\frac{q(3+\eta)}{q+2}}+\Big(n\sum_{i=1}^{n} \big|E Z_1 Z_i\big|\Big)^{\frac{q}{2}}\Bigg\}.
	\end{eqnarray*}
\end{lemma}

\begin{lemma}[ \citet{newman1981}, (12)] \label{lemA3} 	Suppose that $X_1,\ldots, X_m$ are associated, mean zero, finite variance, random variables. Then for any real number $\lambda >0$
	\begin{eqnarray*}
		P\Big(\max\big\{|S_1|,\ldots,|S_m|\big\}\geq   \lambda s_m\Big)\leq 	2 P\Big(|S_m|\geq (\lambda - \sqrt{2}) s_m\Big),
	\end{eqnarray*}
	where $S_m=\sum_{i=1}^{m} X_i$ and $s_m^2= E S_m^2$.
\end{lemma}

\begin{lemma}[ \citet{birkel1988}, Theorem 2]\label{lemA2}
	Let $(X_i)_{i\in\mathbb{N}}$ be a sequence of random variables satisfying $E X_i=0$ and $|X_i|\leq C <\infty$ for $i\in \mathbb{N}$. Assume for some $r>2$
	\begin{eqnarray*}\
		\sup_{k\in \mathbb{N} } \sum_{i : |i-k|\geq n}\mathrm{Cov}(X_{i}, X_{k})= O(n^{-\frac{r-2}{2}})\quad n\in \mathbb{N}.
	\end{eqnarray*}
	Then there is a constant $B$ not depending on $n$ such that for all $n \geq 1$
	\begin{eqnarray*}
		\sup_{m\in \mathbb{N} } E	\left| \sum_{i=m+1}^{m +n } X_i \right|^r  \leq B n^{\frac{r}{2}}.
	\end{eqnarray*}
\end{lemma}

\begin{lemma}[ \citet{moricz1981}, Corollary 1] \label{lemA1} 	Let $\alpha>1$, $\gamma\geq 1$ and $d\geq 1$. Let $\{\xi_i, \, i\in \mathbb{N}^{d}\}$ be real random fields having finite moments of order $\gamma$. Suppose that there exists a nonnegative  and superadditive function $g(R_{b,p})$ of  the rectangle 
		$$R_{b,p}=\big\{(i_1,\ldots,i_d)\in \mathbb{N}^{d}: b_j+1\leq i_j \leq b_j+p_j, \, j=1,\ldots, d\big\}, $$  where $b_j\in \mathbb{N}$ and $p_j\geq 1$, $j=1,\ldots, d$, such that for every $R_{b,p}$ we have
		\begin{eqnarray*}
			E\Big|\sum_{i\in R_{b,p}} \xi_i\Big|^{\gamma}\leq \big(g(R_{b,p})\big)^{\alpha}.
		\end{eqnarray*}
		Then for every $R_{b,p}$ we have
		\begin{eqnarray*}
			E\Big(\max_{1\leq p_1 \leq m_1}\ldots\max_{1\leq p_d\leq m_d}\Big|\sum_{i_1=b_1+1}^{b_1+p_1}\ldots \sum_{i_d=b_d+1}^{b_d+p_d}\xi_{i_1,\ldots,i_d}\Big|\Big)^{\gamma}\leq C(\alpha,\gamma,d) \big(g(R_{b,m})\big)^{\alpha},
		\end{eqnarray*}
		where $C(\alpha,\gamma,d)=(\frac{5}{2})^d \big(1-2^{(1-\alpha)\gamma}\big)^{-d\gamma}$.
	\end{lemma}
	A detailed proof of this lemma is given in \citet{bulinski2007}.


\bibliographystyle{natbib}
\bibliography{asy}

\end{document}